\documentclass[10pt]{article}
\usepackage[dvips]{graphicx}
\usepackage{amssymb}
   \newtheorem{pr}{Proposition}[section]
   \newtheorem{thm}[pr]{Theorem}
   \newtheorem{lemm}[pr]{Lemma}
   
   \newcommand{\qed}{$\rule{2mm}{2mm}$\par\medskip}
   \newcommand{\ff}{\bar f}
   \newcommand{\hh}{\bar h}
   \newcommand{\Rn}{\mathbb{R}}
   \newcommand{\Zn}{\mathbb{Z}}
   \newcommand{\Cn}{\mathbb{C}}

\title{Geometry of Pinchuk's map
}
\author{Janusz Gwo\'zdziewicz}
\date{February 25, 1999}

\begin{document}
\maketitle

\section*{Abstract}
  Sergey Pinchuk found a polynomial map from the real plane to
  itself which is a local diffeomorphism but is not one-to-one.
  The aim of this paper is to give a geometric description of
  Pinchuk's map.

\section{Introduction}
  In the paper~\cite{P} Pinchuk gave an example of a polynomial  map
  $F:\Rn^2\to\Rn^2$ with a non-vanishing Jacobian such that $F$ is not a
  global diffeomorphism.
  Let us recall his construction.  We define the following auxiliary
  polynomials in variables $x$, $y$:
\begin{equation}\label{E1}
  t=xy-1,\quad   h=t(xt+1), \quad f=(xt+1)^2(t^2+y)
\end{equation}
  and a polynomial 
$$ u(f,h)=(1/4)f(75f^3+300f^2h+450fh^2+276f^2+828fh+48h^2+364f+48h). $$
  Then $F=(p,q)$ where $p$ and $q$ are given by
 \begin{eqnarray}\label{E2.1}
    p &=& f+h, \\
    q &=& -t^2-6th(h+1)+u(f,h). \label{E2.2}
 \end{eqnarray}

  Our aim is to give a geometric description of Pinchuk's map. As
  in~\cite{Ca} we find the set of points at which $F$ is not proper.
  It allows us to divide the image of~$F$ into sets with constant multiplicity of
  fibers.  Finally we illustrate how $F$ transforms the real plane to these
  sets.

\section{Geometry of Pinchuk's map}

  Let $f:X\rightarrow Y$ be a continuous map of locally compact spaces. We say
  that the mapping $f$ is not proper at a point $y\in Y$, if there is no a
  neighborhood $U$ of a point $y$ such that the set $f^{-1}(cl(U))$ is compact.

  The set $S_f$ of points at which the map $f$ is not proper
  indicates how the map $f$ differs from a proper map. In particular the
  restriction of $f$ from $X\setminus f^{-1}(S_f)$ to $Y\setminus S_f$ is
  proper. This notion was studied in \cite{jel1} in the complex case and
  in \cite{Pe}, where $S_f$ is called an asymptotic variety.

  Let $C=\Phi(\Rn)$ be a curve given by
$$   \Phi:\Rn\ni s\to(s^2+2s,u(s^2+s,s))\in\Rn^2
$$

\begin{thm}\label{P1}
  The set of points at which $F$ is not proper is the curve~$C$.  The set
  $F^{-1}(C)$ is a smooth curve having three connected components.
  Furthermore
\begin{itemize}
\item[(i)]  for every $v\in\Rn^2\setminus C$ \quad $\#F^{-1}(v)=2$,
\item[(ii)] for every $v\in C\setminus\{(-1,0),(0,0)\}$ \quad $\#F^{-1}(v)=1$,
\item[(iii)] $F^{-1}(-1,0)=F^{-1}(0,0)=\emptyset$.
\end{itemize}
\end{thm}

  A behavior of $F$ is illustrated in figure below.  A set $\Rn^2\setminus
  F^{-1}(C)$ has 4 connected components.  The marked ones are mapped
  diffeomorphically to a marked area above the curve $C$.  Two unmarked regions
  are mapped diffeomorphically to an area below~$C$. Indicated connected
  components of a curve $F^{-1}(C)$ are mapped to indicated parts of the curve
  $C$, respectively.

\noindent
\includegraphics[angle=-90, scale=0.25]{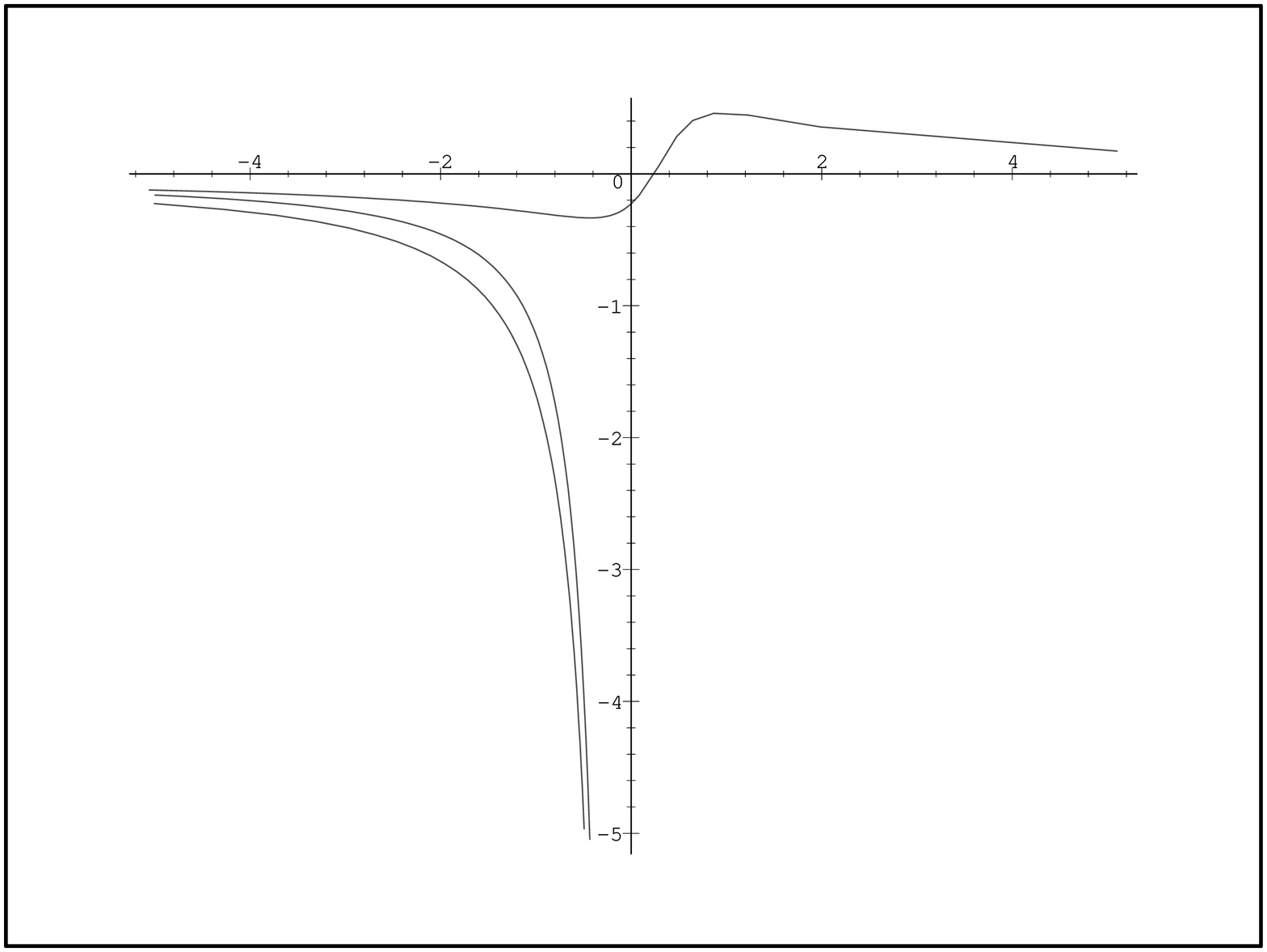}
\includegraphics[angle=-90, scale=0.25]{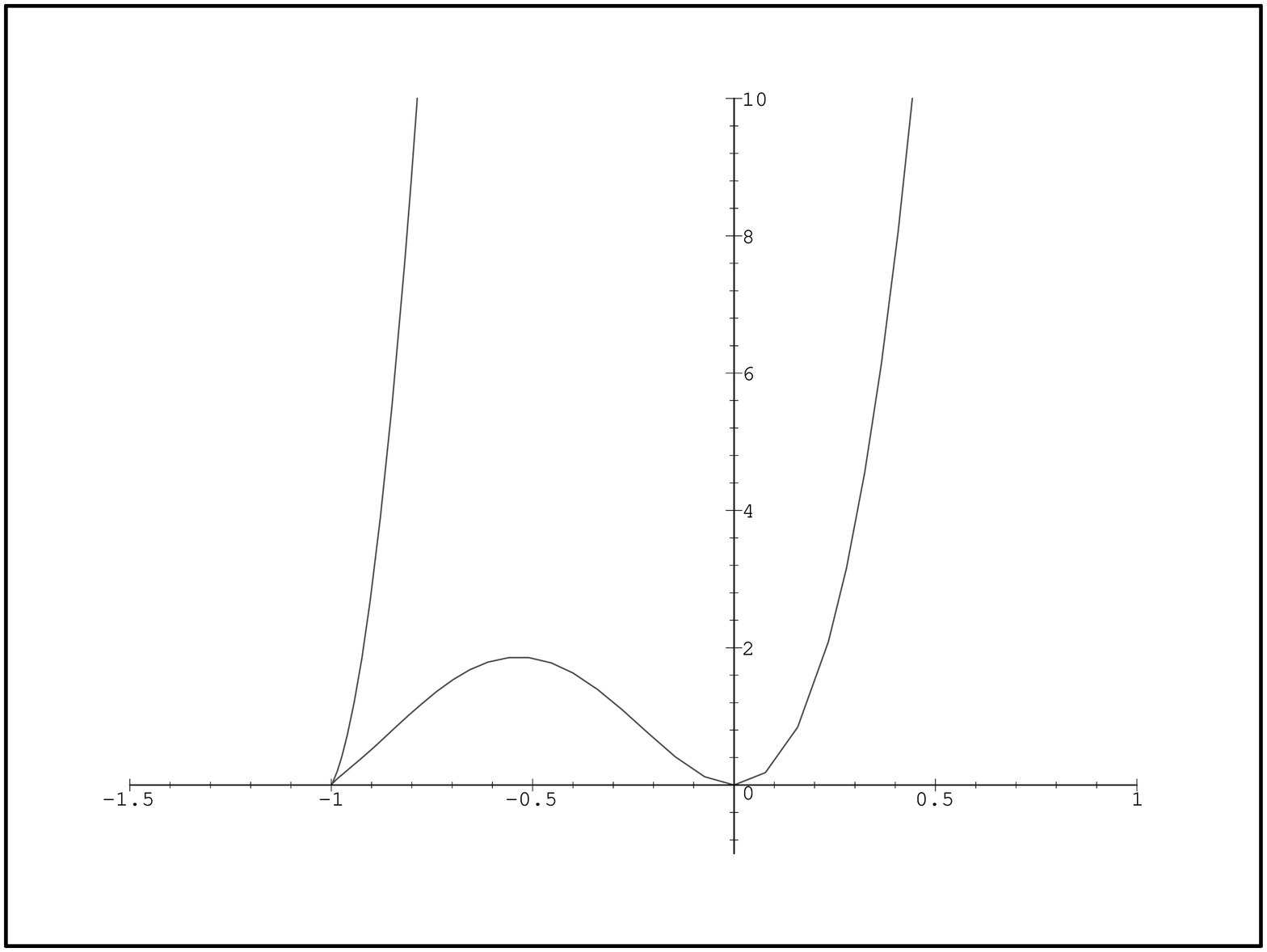}

  A similar description of Pinchuk's map was given in~\cite{Ca}. The author
  characterizes there another mapping from the class discovered by Pinchuk.

  Note that $C$ is not an algebraic curve although it has a polynomial
  pa\-ra\-me\-teri\-za\-tion. It follows from the proof of Lemma~\ref{L4} that the
  Zariski closure of $C$ in $\Rn^2$ consists of~$C$ and a point~$\Phi(s)$,
  where~$s$ is a complex solution of an equation $75s^2+150s+104=0$.

  This shows that Theorem~3 in~\cite{Pe}, where it is claimed that the asymptotic
  variety is algebraic, is partially false.  A polynomial map
  $(x,y)\to(x^2,xy)$ is another counterexample.  An asymptotic variety of this
  map is a right half line $[0,\infty)\times\{0\}$.  Moreover points (b) and
  (c) in Theorem~3 can be removed.  We refer the interested reader
  to~\cite{GJ}, Theorem~3.3.

  \newpage Our approach is based on the following observation due to Zbigniew Jelonek.

\begin{pr}\label{J1}
  Let $\Cn F:\Cn^2\to\Cn^2$ be a complexification of $F$.  Then the set $S_F$
  of points at which $F$ is not proper is contained in the set $S_{\Cn F}$ of
  points at which the mapping $\Cn F$ is not proper.
\end{pr}

\textbf{Proof.} Easy exercise. \qed

The next result concerning complex polynomial maps was proved in~\cite{jel},
\cite{jel1}:

\begin{thm}\label{J2}
  Let $f:\Cn^2\to\Cn^2$ be a dominant polynomial map. Then
  the set $S_f$ of points at which $f$ is not proper consists of a
  finite number (possibly 0) of affine algebraic curves.
\end{thm}

To prove Theorem~\ref{P1} we find the set of points at which complex
Pinchuk's map is not proper and then apply Proposition~\ref{J1}.

\begin{pr}\label{LC}
   $S_{\Cn F}=\Phi(\Cn)$.
\end{pr}

\section{Proofs}

\begin{lemm}\label{L2}
  The following relations hold for $t$, $h$, $f$, $q$:
\begin{eqnarray}
    (h-t)f &=& h^2(h+1)  \label{E5} \\
  Q(f,h,q) &=& 0         \label{E6}
\end{eqnarray}
  where $Q(f,h,q)=f^2(q-u(f,h))+h^2(f-h(h+1))(f+(6f-h)(h+1))$.
\end{lemm}

\textbf{Proof.} The formula~(\ref{E5}) follows immediately from~(\ref{E1}).
  Multiplying equation~(\ref{E2.2}) by $f^2$ gives
  $f^2(q-u(f,h))+(tf)^2+6(tf)h(h+1)f=0$.  By~(\ref{E5}) we have
  $tf=h(f-h(h+1))$. Substituting the right--hand side into above formula we
  obtain (\ref{E6}).
\qed

  Consider a system of equations
\begin{equation}\label{E7}
\left\{ \begin{array}{c}
                              \ff+\hh = a        \\
                         Q(\ff,\hh,b) = 0
\end{array} \right.
\end{equation}

  To solve this system for a given $(a,b)$ it suffices to find roots
  of $W(\ff,a,b)$ where $W(\ff,a,b)=Q(\ff,a-\ff,b)$ is a monic polynomial of
  degree~6 with respect to $\ff$. Therefore a number of solutions
  of~(\ref{E7}) equals a number of roots of a polynomial $W(\ff,a,b)$.

  For every $(\ff,\hh)$ with $\ff\neq0$ there is a unique pair $(a,b)$ such
  that $(\ff,\hh,a,b)$ satisfies~(\ref{E7}).  Namely
\begin{equation}\label{E72}
\left\{ \begin{array}{l}
         a = \ff+\hh        \\
         b = u(\ff,\hh)-\hh^2(\ff-\hh(\hh+1))(\ff+(6\ff-\hh)(\hh+1))/\ff^2
\end{array} \right.
\end{equation}
  Observe that, generically, the Pinchuk
  mapping is a composition of a polynomial map  $(f,h)$ and a map $G$ given
  by~(\ref{E72}). We precise this idea in Lemmas~\ref{L5} and~\ref{L51}. First
  we need some preparation.

  \newpage
  Let $B=\{\,(\ff,\hh)\in\Cn^2: \ff(\ff-\hh(\hh+1))=0\,\}$.

\begin{lemm}\label{L3}
  The restriction of $(f,h):\Cn^2\to\Cn^2$ from $\Cn^2\setminus
  f^{-1}(0)$ to $\Cn^2\setminus B$ is a homeomorphism.  A curve $f^{-1}(0)$
  has two algebraic components:  
  \newline$A_1=(f,h)^{-1}(0,0)$ and $A_2=(f,h)^{-1}(0,-1)$.
  The restriction of the polynomial function~$t$ to $A_i$ is injective and 
  $t(A_i)=\Cn\setminus\{0\}$ for $i=1,2$.
\end{lemm}

\textbf{Proof.}
  We claim that $(f,h)^{-1}(B)=f^{-1}(0)$. Indeed, suppose that
  $f\neq0$ and $f-h(h+1)=0$.  Substituting  $f=h(h+1)$ into
  (\ref{E5}) gives $ft=0$.  Hence $t=0$ and by~(\ref{E1})  $h=0$.
  Therefore $f=h(h+1)=0$ -- contrary to our assumption.

   An easy computation shows that the mapping
   $\Cn^2\setminus B\ni(\ff,\hh)\to((\hh+1)\ff/(\ff-\hh(\hh+1))^2,
   (\ff-\hh^2)(\ff-\hh(\hh+1))^2/\ff^2)\in\Cn^2$ composed with $(f,g)$ is an
   identity on $\Cn^2\setminus B$. This together with the claim gives the
   first part of the lemma.

  Now we describe a curve $f^{-1}(0)$. Assume that $f(x,y)=0$.
  Then by~(\ref{E1}) $xt+1=0$ or $t^2+y=0$.  We obtain by~(\ref{E1}) in
  the first case $x=-1/t$, $y=-t(t+1)$, $h=0$ and in the second case
  $x=-(t+1)/t^2$, $y=-t^2$, $h=-1$.  This finishes the proof.
\qed

\begin{lemm}\label{L5}
  The complex Pinchuk mapping $\Cn F$ restricted to $f^{-1}(0)$
  is two-to-one between $f^{-1}(0)$ and $\{-1,0\}\times(\Cn\setminus\{0\})$.
\end{lemm}

\textbf{Proof.} By (\ref{L3})
  $q(x,y)=-t^2-6th(h+1)+u(f,h)=-t^2$ for $(x,y)\in f^{-1}(0)$.  Thus by the
  second part of Lemma~\ref{L3} mappings
  $(f,h)^{-1}(0,0)\to\{0\}\times(\Cn\setminus\{0\})$ and
  $(f,h)^{-1}(0,-1)\to\{-1\}\times(\Cn\setminus\{0\})$ induced by $\Cn F$ are
  two-to-one.
\qed

\begin{lemm}\label{L51}
  The complex Pinchuk mapping $\Cn F$ restricted to $\Cn^2\setminus f^{-1}(0)$
  is a composition
$$ \Cn^2\setminus f^{-1}(0) \stackrel{(f,h)}{\longrightarrow} \Cn^2\setminus B
   \stackrel{G}{\longrightarrow} \Cn^2
$$
  where $(a,b)=G(\ff,\hh)$ is given by~(\ref{E72}).
\end{lemm}

\textbf{Proof.} It follows immediately from \ref{L2} and \ref{L3}.
\qed

\begin{lemm}\label{L6}
  Let $(\ff,\hh,a,b)$ satisfy~(\ref{E7}) and let $(\ff,\hh)\in B$.  Then
  \begin{itemize}
    \item[(i)]  if $\ff=0$ then $a=-1$ or $a=0$,
    \item[(ii)] if $\ff\neq0$, $\ff=\hh(\hh+1)$ then $(a,b)=\Phi(\hh)$.
  \end{itemize}
\end{lemm}

\textbf{Proof.}
  Putting $\ff=0$ in~(\ref{E7}) we get $\hh=a$, $Q(0,\hh,b)=\hh^4(\hh+1)^2$
  and~(i) follows.  To prove~(ii) put $\ff=\hh(\hh+1)$ in~(\ref{E72}).  We
  obtain $a=\ff+\hh$, $b=u(\ff,\hh)$.  Thus $(a,b)=\Phi(\hh)$.
\qed

\textbf{Proof of Proposition~\ref{LC}.}  First we show that  
$$ S_{\Cn F}\subset(\{-1\}\times\Cn)\cup (\{0\}\times\Cn)\cup \Phi(\Cn). $$

  Fix
  $(a_0,b_0)\in\Cn^2$ lying in the complement of above curves and take a
  neighborhood $U$ of $(a_0,b_0)$ such that $cl(U)$ is a compact set disjoint
  from $(\{-1\}\times\Cn) \cup(\{0\}\times\Cn)\cup \Phi(\Cn)$. Let $V$ be the
  set of all $(\ff,\hh)\in\Cn^2$ such that $(\ff,\hh,a,b)$
  satisfies~(\ref{E7}) for some $(a,b)\in cl(U)$.  By continuity of roots (see
  \cite{rie}, Proposition~1.5.5) $V$ is a compact set.  Moreover, it follows
  from Lemma~\ref{L6} that $V\cap B=\emptyset$.  Hence $V=G^{-1}(cl(U))$.
  Since $\Cn F^{-1}(cl(U))=(f,h)^{-1}(V)$, we conclude by Lemma~\ref{L3} that
  $\Cn F^{-1}(cl(U))$ is compact.  Therefore the map $\Cn F$ is proper at the
  point $(a_0,b_0)$.

  Now we show that the polynomial $W(\ff,a,b)$ has no multiple factor in a
  ring $\Cn[\ff,a,b]$. Assume that it has. Then
  $Q(\ff,\hh,b)=W(\ff,\ff+\hh,b)$ has also a multiple factor. Since  $Q$ can
  be rewritten as $Q(\ff,\hh,b)=\ff^2b+Q_1(\ff,\hh)$ the only candidate for a
  multiple factor is $\ff$.  But $\ff$ does not divide $Q(\ff,\hh,b)$ contrary
  to our assumption.

  Thus for a generic~$(a,b)\in\Cn^2$ the polynomial $W(\ff,a,b)$ has 6 single
  roots.  It follows that $\#\Cn F^{-1}(a,b)=6$ for $(a,b)$ in a general
  position.

  Our next task is to show that lines $\{0\}\times\Cn$ and
  $\{-1\}\times\Cn$ are not contained in $S_{\Cn F}$.  We will use the
  following property (see \cite{jel1}):  \textit{If $\#\Cn F^{-1}(a,b)=\#\Cn
  F^{-1}(a_{\mathrm{gen}},b_{\mathrm{gen})}$, then $\Cn F$ is proper at a
  point~$(a,b)$.}  Here $(a_{\mathrm{gen}},b_{\mathrm{gen})}$ denotes a point
  in a general position, so $\#\Cn
  F^{-1}(a_{\mathrm{gen}},b_{\mathrm{gen})}=6$.

  We have $W(\ff,0,b)=\ff^2(-197/4\ff^4+104\ff^3-63\ff^2+b)$.  Thus for a
  generic $b\in\Cn$ a system~(\ref{E7}) has 4 solutions $(\ff_i,\hh_i)$
  ($i=1,\dots,4$) such that $(\ff_i,\hh_i)\notin B$ and one double solution
  $(\ff_5,\hh_5)=(0,0)$.  By~\ref{L5} and~\ref{L51} we get
  $\#\Cn F^{-1}(0,b)=6$.  Hence by the above mentioned property
  $\{0\}\times\Cn\not\subset S_{\Cn F}$.   We show in the same way that
  $\{-1\}\times\Cn\not\subset S_{\Cn F}$. From \ref{J2}
  it follows that the set $S_{\Cn F}$ is either a curve
  $\Phi(\Cn)$ or is empty.  The latter is impossible because in this case
   real Pinchuk's map would be proper and thus it would be a diffeomorphism.
\qed

  From now we treat all polynomials under considerations as real polynomials.

\begin{lemm}\label{L4}
   The parameterization $\Phi$ of a curve $C$ is injective and
   $C$ is smooth except a point $(-1,0)$.
\end{lemm}

\textbf{Proof.} Suppose to the contrary that there are $s_1\neq s_2$ such that
  $\Phi(s_1)=\Phi(s_2)$.  An equation $s_1^2+2s_1=s_2^2+2s_2$ gives
  $s_2=-2-s_1$. Combining this with equality  $u(s_1^2+s_1,s_1)=
  u(s_2^2+s_2,s_2)$ we get $(75s_1^2+150s_1+104)(s_1+1)^3=0$.  Hence
  $s_1=s_2=-1$ which contradicts the choice of $s_1$ and $s_2$.

  Since $\frac{d}{ds}\Phi(s)$ vanishes for $s=-1$ only,
  the curve $C$ has a singularity at most at~$(-1,0)$.
\qed

\begin{lemm}\label{L1}
  Let $H:\Rn^n\to\Rn^n$ be a  polynomial function which is a local diffeomorphism.
  Then the function  $\Rn^n\ni v\to\#H^{-1}(v)\in\Zn$ is lower semi-continuous.
  Moreover, $H$ is proper at a point $v_0\in\Rn^n$ if and only if 
  $\#H^{-1}(v)$ is locally constant at $v_0$.
\end{lemm}

   The purely topological proof is left to the reader.

\smallskip
\textbf{Proof of Theorem~\ref{P1}.}
   It follows from~\ref{J1} and~\ref{LC} that $S_F$ is a subset of
   $\Phi(\Cn)\cap\Rn^2$. Since an asymptotic variety does not have isolated
   points (see~\cite{Pe}, Theorem~3), $S_F=\Phi(\Rn)=C$.  From the topological
   point of view $C\subset\Rn^2$ is an embedded line.  Hence $\Rn^2\setminus
   C$ has 2 connected components.  We shall calculate a multiplicity of fibers
   in each of them.

   The curve $C$ cuts a line $\{3\}\times\Rn$ transversely at points
   $(3,3142)=\Phi(1)$ and $(3,8406)=\Phi(-3)$.  Hence the points $(3,0)$ and
   $(3,4000)$ lie on opposite sides of $C$.  One checks that polynomials 
   $W(\ff,3,0)$ and $W(\ff,3,4000)$ have 2 real roots (e.g. we can use Sturm sequence).  Hence
   $\#F^{-1}(3,0)=\#F^{-1}(3,4000)=2$ and by Lemma~\ref{L1} $\#F^{-1}(v)=2$
   for every $v\in\Rn^2\setminus C$.

   Now we show that $F^{-1}(0,0)=F^{-1}(-1,0)=\emptyset$.  We have
   $W(\ff,0,0)=-\ff^4(197/4\ff^2-104\ff+63\ff)$.  The only real root of this
   polynomial is $\ff=0$.  Since by \ref{L5} $(0,0)$ does not belong to the
   set $F(f^{-1}(0))$ we have $F^{-1}(0,0)=\emptyset$.  We check similarly
   that $F^{-1}(-1,0)=\emptyset$.

   It remains to compute $\#F^{-1}(v)$ for $v\in C$.  Let $D(a,b)$ be a
   discriminant of the polynomial $W(\ff,a,b)$ with respect to $\ff$.  Using
   any  computer symbolic algebra program one can check that the polynomial
   $D$ is nonzero at a point $\Phi(1)\in C$. Hence $D$ does not vanish on $C$
   but a finite number of points.  Fix $(a_0,b_0)=\Phi(s)\in C$ for which
   $D(a_0,b_0)\neq0$.  We have shown before that $W(\ff,a,b)$ has two real
   roots for $(a,b)\in\Rn^2\setminus C$.  Therefore by continuity of roots the
   system~(\ref{E7}) has two real solutions at a point $(a_0,b_0)$.  By
   \ref{L6} and \ref{L4} one of these solutions is $(\ff,\hh)=(s(s+1),s)\in
   B$.  Thus $\#F^{-1}(a_0,b_0)=\#(G^{-1}(a_0,b_0)\cap\Rn^2)=1$.

  By Lemma~\ref{L1} we see that $\#F^{-1}(v)\leq1$ for $v\in C$. Let
  $$ K=\{v\in C:\#F^{-1}(v)=0\}. $$  We have shown that $K$ is a finite set and
  points $(-1,0)$, $(0,0)$ belong to $K$.

  The restriction of $F$ from $F^{-1}(C)$ to $C\setminus K$ is a
  diffeomorphism.  Hence the curve $F^{-1}(C)$ has $\#K+1$ connected
  components homeomorphic to a line and  its complement $\Rn^2\setminus
  F^{-1}(C)$ has $\#K+2$ connected components.

  The set $\Rn^2\setminus C$ has two connected components $S_1$, $S_2$. Write
  $F^{-1}(S_i)=\bigcup_{j=1}^{r_i}S_{i,j}$ where $S_{i,j}$ are connected
  components of the set $F^{-1}(S_i)$.  The mappings $S_{i,j}\to S_i$
  ($i=1,2$, $j=1,\dots,r_i$) induced by $F$ are proper and unramified hence
  they are topological coverings.  Since sets $S_i$ ($i=1,2$) are simply
  connected, these mappings are homeomorphisms.  Since $\#F^{-1}(v)=2$ for
  $v\in\Rn^2\setminus C$,\  $F^{-1}(S_i)$ consists of 2 connected components
  for $i=1,2$.  It follows that $\#K=2$, therefore $K=\{(-1,0),(0,0)\}$.
  \qed

 \end{document}